\documentclass[9pt]{amsart}
\usepackage{amsfonts}
\usepackage{graphicx}
\usepackage[all]{xy}
\usepackage{amssymb,extarrows}
\usepackage{amscd}
\usepackage{amsmath,amsfonts,amssymb,amsthm,mathrsfs,xypic}
\newtheorem{thm}{Theorem}[section]

\newtheorem{cor}[thm]{Corollary}
\newtheorem{lem}[thm]{Lemma}
\newtheorem{exm}[thm]{Example}

\newtheorem{prop}[thm]{Proposition}

\newtheorem{defn}[thm]{Definition}
\numberwithin{equation}{section}

\newcommand{\s}{\hfill\blacksquare}
\newcommand{\End}{\operatorname{End}}
\newcommand{\add}{\operatorname{add}}
\newcommand{\Ima}{\operatorname{Im}}
\newcommand{\Hom}{\operatorname{Hom}}
\newcommand{\Ext}{\operatorname{Ext}}

\begin{document}
\title [Categorical resolution]{Categorical resolutions \\ of a class of derived categories}
\author [Pu Zhang]
{Pu Zhang}
\thanks{Supported by NSF China (Grant No. 11271251
and 11431010).}
\thanks {pzhang$\symbol{64}$sjtu.edu.cn}
\address{Department of Mathematics \\ Shanghai Jiao Tong
University, Shanghai 200240, China}
\begin{abstract} \ Using the relative derived categories,
we prove that if an Artin algebra $A$ has a module $T$ with ${\rm
inj.dim}T<\infty$ such that $^\perp T$ is finite, then the bounded
derived category $D^b({\rm mod}A)$ admits a categorical resolution;
and that for CM-finite Gorenstein algebra, such a categorical
resolution is weakly crepant.

\vskip5pt


\vskip5pt

\end{abstract}
\maketitle
\section {\bf Introduction and preliminaries}
\subsection{} A resolution of a
singular variety motivates the study of a categorical resolution of
a triangulated category $\mathcal D$, which is a pair $(D^b(\mathcal
A), \mathcal K)$, where $\mathcal A$ is an abelian category of
finite homological dimension, $\mathcal K$ a thick subcategory of
$D^b(\mathcal A)$ such that $D\cong D^b(\mathcal A)/\mathcal K$ (A.
Bondal and D. Orlov [BO]). In [K] A. Kuznetsov defines a categorical
resolution of $\mathcal D$ as a triple $(\widetilde{\mathcal D},
\pi_*, \pi^*),$ where $\widetilde{\mathcal D}$ is an admissible
subcategory of $D^b(\widetilde{X})$ with $\widetilde{X}$ a smooth
variety, $\pi_*: \widetilde{\mathcal D} \rightarrow \mathcal D$ and
$\pi^*: \mathcal D_{\rm perf} \rightarrow \widetilde{\mathcal D}$
are triangle functors satisfying some conditions. A non-commutative
crepant resolution in the sense of M. Van den Bergh [VB] has all
these properties. For other influential works see e.g. [BKR], [BLV],
[Lu] and [SV].

\vskip5pt

In this paper we combine Kuznetsov's definition with Bondal-Orlov's
one (Definition \ref{defcatres}). Throughout $\mathcal A$ is an
abelian category and $A$ an Artin algebra. If $\mathcal A$ has
enough projective objects, then $\mathcal P: = \mathcal P(\mathcal
A)$ denotes the full subcategory of projective objects.

\vskip5pt

For a full additive subcategory $\mathcal C$ of $\mathcal A$, the
relative derived category $D_\mathcal C^b(\mathcal A)$ has been
studied in different setting up (e.g. [N], [Bu], [GZ], [C2] and
[AHV]). If $\mathcal C$ is a resolving and contravariantly finite,
then $D^b(\mathcal A)$ can be described as $D_\mathcal C^b(\mathcal
A)/K^b_{ac}(\mathcal C)$ (Proposition \ref{relativedescription}). So
we get the Verdier functor $\pi_*: D_\mathcal C^b(\mathcal
A)\rightarrow D^b(\mathcal A)$. If $D^b_{\rm perf}(\mathcal A) =
K^b(\mathcal P)$, then we get a triple $(D_\mathcal C^b(\mathcal A),
\pi_*, \pi^*)$ with embedding $\pi^*: D^b_{\rm perf}(\mathcal
A)\rightarrow D_\mathcal C^b(\mathcal A)$. This will be served as a
categorical resolution of $D^b(\mathcal A)$ in our consideration.

\subsection{} Let ${\rm Mod}A$ (resp. ${\rm mod}A$) be the
category of right (resp. finitely generated) $A$-modules. For $T\in
{\rm mod}A$, let $^\perp T$ be the full subcategory of ${\rm mod}A$
consisting of $A$-modules $X$ with $\Ext^i_A(X, T) = 0, \ \forall \
i\ge 1;$ and ${\rm add} T$ (resp. ${\rm Add} T$) the full
subcategory of {\rm mod}$A$ (resp. {\rm Mod}$A$) consisting of
direct summands of finite (resp. arbitrary) direct sums of copies of
$T$. Let \ $^{\perp_{\rm big}} ({\rm Add} T)$ the full subcategory
of $A$-Mod given by $\{X\in {\rm Mod}A \ | \ \Ext^i_A(X, T') = 0, \
\forall \ i\ge 1, \ \forall \ T'\in {\rm Add} T \}$.

\vskip5pt

Assume that ${\rm gl.dim} A=\infty$, and there are modules $T$ and
$M$ in mod$A$ with ${\rm inj.dim} T< \infty$ such that $^\perp T =
{\rm add} M$.  The main result Theorem \ref{mainresult} claims that
the above triple $(D^b({\rm mod}B), \pi_*, \pi^*)$ is a categorical
resolutions of $D^b({\rm mod}A)$; and that if $A$ is a {\rm
CM}-finite Gorenstein algebra and $B$ its relative Auslander
algebra, then $D^b({\rm mod}B)$ is a weakly crepant categorical
resolution of $D^b({\rm mod}A)$. The same result holds also for
$D^b({\rm Mod}A)$.

\subsection{} An object $P$ of triangulated
category $\mathcal D$ is {\it perfect} ([O1], [K]), if for each
$Y\in\mathcal D$ there are only finitely many $i\in \Bbb Z$ with
$\Hom(P, Y[i]) \ne 0$. Let $\mathcal D_{\rm perf}$ be the full
subcategory of perfect objects. It is a thick subcategory of
$\mathcal D$, and an invariant of triangle-equivalences. If
$\mathcal A$ has enough projective objects, then $K^b(\mathcal P)
\subseteq D^b_{\rm perf}(\mathcal A)$; but we stress that $D^b_{\rm
perf}(\mathcal A) = K^b(\mathcal P)$ is not always true, although
this is the case in the many situations (Section 2).

\vskip5pt

Define {\it the singularity category} of $\mathcal A$ to be
$D^b_{sg}(\mathcal A): = D^b(\mathcal A)/D^b_{\rm perf}(\mathcal A)$
([B], [O2]). A triangulated category $\mathcal D$ is {\it smooth}
([BO], [K]), if it is triangle-equivalent to $D^b(\mathcal A)$ with
$D^b_{\rm sg}(\mathcal A) = 0$. For other definitions of the
smothness see e.g. [KS], [Lu] and [TV].

\begin{defn} \label{defcatres}  $(1)$ \ {\rm ([BO],  [K, 3.2])} \ A categorical resolution of a non-smooth triangulated
category $\mathcal D$ is a triple $(\widetilde{\mathcal D}, \pi_*,
\pi^*)$,  where  $\widetilde{\mathcal D}$ is a smooth triangulated
category,  $\pi_*: \widetilde{\mathcal D} \rightarrow \mathcal D$
and $\pi^*: \mathcal D_{\rm perf} \rightarrow \widetilde{\mathcal
D}$ are triangle functors, such that

${\rm(i)}$ \ $\pi_*$ induces a triangle-equivalence
$\widetilde{\mathcal D}/{\rm Ker}\pi_*\cong \mathcal D;$

${\rm(ii)}$ \ $\pi^*$ is left adjoint to $\pi_*$ on $\mathcal D_{\rm
perf}$, that is, there is a functorial isomorphism $\eta_{P, X}:
\Hom_{\widetilde{\mathcal D}}(\pi^*P, X)\cong \Hom_\mathcal D(P,
\pi_*X), \ \ \forall \ P\in \mathcal D_{\rm perf}, \ \forall \ X\in
\widetilde{\mathcal D};$

${\rm(iii)}$ \ The unit $\eta = (\eta_P)_{P\in \mathcal D_{\rm
perf}}: {\rm Id}_{\mathcal D_{\rm perf}} \rightarrow \pi_*\pi^*$ is
a natural isomorphism of functors, where $\eta_P$ is the morphism
$\eta_{P, \pi^*P}({\rm Id}_{\pi^*P}): P\rightarrow \pi_*\pi^*P$.

\vskip5pt

 $(2)$ \ {\rm ([K, 3.4])} \ A categorical resolution $(\widetilde{\mathcal D}, \pi_*,
\pi^*)$ of a triangulated category $\mathcal D$ is weakly crepant if
$\pi^*$ is also right adjoint to $\pi_*$ on $\mathcal D_{\rm perf}$.
\end{defn}

Note that ${\rm(iii)}$ implies that $\pi^*: \mathcal D_{\rm perf}
\rightarrow \widetilde{\mathcal D}$ is fully faithful. If $\pi_*:
\widetilde{\mathcal D} \rightarrow \mathcal D$ is full and dense,
then ${\rm(i)}$ holds automatically; however $\pi_*$ usually can not
be full.

\vskip5pt

It is well-known that for a complex singular variety $X$ there is a
proper birational resolution of singularities
$\widetilde{X}\rightarrow X$; and if $X$ is of rational singularity,
then $D^b(\widetilde{X})$ is a categorical resolution of $D^b(X)$. A
non-commutative crepant resolution ([VB]) induces a weakly crepant
categorical resolution ([K]).

\subsection{} Let $\mathcal B$ be an additive category,  $\mathcal C$ its additive full
subcategory, and $X\in \mathcal B$. A morphism $f: C\rightarrow X$
with $C\in\mathcal C$  is {\it a right $\mathcal C$-approximation}
of $X$, if the induced map $\Hom_\mathcal B(C', f)$ is surjective
for each $C'\in\mathcal C$. If each object of $\mathcal B$ has a
right $\mathcal C$-approximation, then $\mathcal C$ is {\it
contravariantly finite} in $\mathcal B$ ([AR]). For example, for
$M\in$ {\rm mod}$A$, ${\rm add} M$ {\rm(}resp. ${\rm Add} M${\rm)}
is contravariantly finite in {\rm mod}$A$ {\rm(}resp. {\rm
Mod}$A${\rm)}.

\vskip5pt

Let $\mathcal A$ be an abelian category with enough projective
objects. A full subcategory $\mathcal C$ is {\it resolving} ([AB]),
provided that $\mathcal C\supseteq \mathcal P,$ $\mathcal C$ is
closed under extensions and direct summands, and that $\mathcal C$
is closed under the kernels of epimorphisms.  A resolving
subcategory is additive.

\vskip5pt

Denote by $\mathcal{GP}(\mathcal A)$ the full subcategory of
$\mathcal A$ consisting of Gorenstein-projective objects ([EJ]).
Then $\mathcal{GP}(\mathcal A)$ is a resolving subcategory of
$\mathcal A$, and $\mathcal{GP}(\mathcal A)$ is contravariantly
finite in $\mathcal A$  ([AR], [EJ], [H]), if each object of
$\mathcal{A}$ has a finite Gorenstein-projective dimension. Note
that  $\mathcal {GP}({\rm Mod}A)$ is contravariantly finite in {\rm
Mod}$A$ {\rm(}{\rm [Be1, 3.5]}{\rm)}. An Artin algebra $A$ is {\it
CM-finite}, if $\mathcal{GP}({\rm mod}A)$ has only finitely many
pairwise non-isomorphic indecomposable objects; and  $A$ is {\it
Gorenstein}, if ${\rm inj.dim} A_A$ $< \infty$ and ${\rm inj.dim}\
_AA < \infty$. For a {\rm CM}-finite algebra $A$, $\mathcal
{GP}({\rm mod}A)$ is contravariantly finite in {\rm mod}$A$. For
examples of \ {\rm CM}-finite non-Gorenstein algebras we refer to
{\rm [Rin]}. If $A$ is virtually Gorenstein, then $\mathcal
{GP}({\rm mod}A)$ is contravariantly finite in $A$-{\rm mod} {\rm
([Be1, 8.2])}. In particular, if $A$ is Gorenstein then $\mathcal
{GP}({\rm mod}A)$ is contravariantly finite in $A$-{\rm mod}.

\section{\bf Perfect subcategory of a triangulated
category}
We give two classes of abelian categories $\mathcal A$
with enough projective objects, such that $D^b_{\rm perf}(\mathcal
A) = K^b(\mathcal P),$ and also an example of $\mathcal A$ such that
$D^b_{\rm perf}(\mathcal A) \ne K^b(\mathcal P).$

\begin{lem}\label{buch} {\rm ([B, 1.2.1])} \ Let $P\in D^b(\mathcal A)$. Then $P\in K^b(\mathcal
P)$ if and only if there is a finite set $I(P)\subseteq \Bbb Z$,
such that $\Hom_{D^b(\mathcal A)}(P, M[j]) = 0$ for $j\notin I(P)$
and each object $M\in\mathcal A$.
\end{lem}
$\mathcal A$ is {\it finitely filtrated}, if there are finitely many
objects $S_1, \cdots, S_m$, such that for $0\ne X\in\mathcal A$
there is a sequence of monomorphisms $0= X_0 \stackrel {f_1}
\longrightarrow \cdots \stackrel {f_n} \longrightarrow X_n = X$ with
each ${\rm Coker} f_i \in\{S_1, \cdots, S_m\}$.  For example, mod$A$
is finitely filtrated for Artin algebra $A$.

\begin{prop} \label{finitely filtrated} Let $\mathcal A$ be an abelian category with
enough projective objects.

$(1)$ \ If $\mathcal A$ is finitely filtrated, then $D^b_{\rm
perf}(\mathcal A) = K^b(\mathcal P)$.

$(2)$ \ If $\mathcal A$ has infinite direct sums, then $D^b_{\rm
perf}(\mathcal A) = K^b(\mathcal P)$.
\end{prop}
\noindent{\bf Proof.} \ We only prove  $D^b_{\rm perf}(\mathcal
A)\subseteq K^b(\mathcal P)$. Let $P\in D^b_{\rm perf}(\mathcal A)$.

\vskip5pt

$(1)$ \ Let $\mathcal A$ be finitely filtrated by $S_1, \cdots,
S_m$. Put $S: = \bigoplus\limits_{1\le i\le n}S_i$. Let $I(P)$ be
the finite set of integers $i$ with $\Hom_{D^b(\mathcal A)}(P, S[i])
\ne 0$. Then $\Hom_{D^b(\mathcal A)}(P, S[j]) = 0, \ \forall \
j\notin I(P)$. Since any object of $\mathcal A$ has a filtration
with factors in $\{S_1, \cdots, S_m\}$, and any short exact sequence
in $\mathcal A$ forms a distinguished triangle in $D^b(\mathcal A)$,
it follows that $\Hom_{D^b(\mathcal A)}(P, M[j]) = 0$ for $j\notin
I(P)$ and any object $M\in\mathcal A$. So  $P\in K^b(\mathcal P)$ by
Lemma \ref{buch}.

\vskip5pt

$(2)$  \  The idea is due to J. Rickard [Ric, Prop. 6.2]. Take a
quasi-isomorphism $Q \rightarrow P$ with $Q\in K^{-, b}(\mathcal P)$
and ${\rm H}^nP = 0, \ \forall \ n\le N$. It suffices to prove that
there is an integer $n$ with $n\le N$ such that ${\rm Im}d^{n}_Q\in
\mathcal P$. Otherwise, ${\rm
 Im}d^n_Q\notin \mathcal P$ for each $n\le N$.  Since $\mathcal A$ has infinite direct
sums, $M:= \bigoplus\limits_{n\le N}{\rm
 Im}d_Q^{n}\in \mathcal A$. The non-zero morphism
$$f: Q^{n}\rightarrow M = \bigoplus\limits_{j\le N}{\rm
 Im}d_Q^{j} = {\rm Im}d_Q^{n}\oplus (\bigoplus\limits_{j\le N, j\ne n}{\rm
 Im}d_Q^{j})$$ induces a chain map $Q
\rightarrow M[-n]$. Since ${\rm Im}d_Q^n\notin \mathcal P$, it
follows that this chain map is not null homotopic. Thus ${\rm
Hom}_{K^-(\mathcal A)}(Q, M[-n])\ne 0,  \ \forall \ n\le N,$  and
hence $${\rm Hom}_{D^b(\mathcal A)}(P, M[-n])\cong {\rm
Hom}_{D^-(\mathcal A)}(Q, M[-n])
 \cong {\rm Hom}_{K^-(\mathcal A)}(Q, M[-n])\ne 0.$$ This contradicts the
assumption $P\in D^b_{\rm perf}(\mathcal A)$. \hfill $\s$

\begin{exm} \label{no=} We include an example such that $K^b(\mathcal P) \subsetneqq D^b_{\rm perf}(\mathcal
A)$. Let $S$ be the polynomial algebra $k[x_0, \cdots, x_n]$ over
field $k$,  and $R$ the exterior algebra $k[x_0, \cdots,
x_n]/\langle x_i^2, x_ix_j+x_jx_i\rangle$. Then $K^b(\mathcal P({\rm
gr}R)) \subsetneqq D^b_{\rm perf}({\rm gr}R)$, where  ${\rm gr}R$ is
the category of finitely generated graded right $A$-modules, which
is an abelian category with enough projective objects.

\vskip5pt

Otherwise $K^b(\mathcal P({\rm gr}R)) = D^b_{\rm perf}({\rm gr}R)$.
Since ${\rm gr}S$ is of finite global dimension, we have
$K^b(\mathcal P({\rm gr}S))  = D^b({\rm gr}S)$, and hence $D^b_{\rm
perf}({\rm gr}S) = D^b({\rm gr}S).$ On the other hand, by the {\rm
BGG} equivalence $D^b({\rm gr}S)\cong D^b({\rm gr}R)$ {\rm([BGG];}
or {\rm [OSS, Thm. 4.5, p.227])} we get $D^b_{\rm perf}({\rm
gr}S)\cong D^b_{\rm perf}({\rm gr}R)$. This deduces $D^b_{\rm
perf}({\rm gr}R) = D^b({\rm gr}R),$ and hence $K^b(\mathcal P({\rm
gr}R)) = D^b({\rm gr}R)$. But this implies that each object in ${\rm
gr}R$ has finite projective dimension, which is absurd.
\end{exm}
\section{\bf Relative derived categories}
\subsection{} Let $\mathcal C$ be a full additive subcategory of abelian category $\mathcal A$. A complex $M^\bullet$ over $\mathcal
A$ is {\it $\mathcal C$-acyclic}, if $\operatorname{Hom}_{\mathcal
A}(C, M^\bullet)$ is acyclic for all objects $C\in \mathcal C$. A
chain map \ $f^\bullet$ is {\it a $\mathcal C$-quasi-isomorphism},
if the induced chain map \ $\operatorname{Hom}_{\mathcal A}(C,
f^\bullet)$ is a quasi-isomorphism for all objects $C\in \mathcal
C$, or equivalently, the mapping cone
$\operatorname{Con}(f^\bullet)$ is $\mathcal C$-acyclic.  For $*\in
\{b, -, {\rm blank}\}$, let $K^*_{\mathcal C ac}(\mathcal A)$ denote
the full subcategory of the homotopy category $K^*(\mathcal A)$
consisting of $\mathcal C$-acyclic complexes. Then
$$K^*_{\mathcal C ac}(\mathcal A) = \ ^\perp\mathcal C: = \{X^\bullet\in K^*(\mathcal A)\ | \
\Hom_{K^*(\mathcal A)}(C, X^\bullet[n]) =0, \ \forall \ n\in\Bbb Z,
\ \forall \ C\in\mathcal C\}.$$ So $K^*_{\mathcal C ac}(\mathcal A)$
is thick in $K^*(\mathcal A)$. {\it The $\mathcal C$-relative
derived category} is the Verdier quotient \ $D^*_{\mathcal
C}(\mathcal A): = K^*(\mathcal A)/K^*_{\mathcal Cac}(\mathcal A)$ \
([AHV], [GZ], [C2]). If $\mathcal A$ has enough projective objects
and $\mathcal C\supseteq \mathcal P$, then $K^*_{\mathcal
Cac}(\mathcal A)$ is a thick subcategory of $K^*_{ac}(\mathcal A)$,
the full subcategory of $K^*(\mathcal A)$ consisting of acyclic
complexes, and hence the derived category $D^*(\mathcal A): =
K^*(\mathcal A)/K^*_{ac}(\mathcal A)$ is the Verdier quotient of
$D^*_{\mathcal C}(\mathcal A)$ by $K^*_{ac}(\mathcal
A)/K^*_{\mathcal Cac}(\mathcal A)$. Thus we have the Verdier functor
$\pi_*: D^*_{\mathcal C}(\mathcal A)\rightarrow D^*(\mathcal A)$,
which is an equivalence if and only if $\mathcal {C} = \mathcal P$.

\vskip5pt

For examples, if $\mathcal A$ has enough projective objects and
$\mathcal C: = \mathcal {GP}(\mathcal A)$, then $D^*_{\mathcal
C}(\mathcal A)$ is the Gorenstein derived category ({\rm  [GZ]}); if
$A$ is an Artin algebra and $M\in A$-mod,  then we have the
$M$-relative derived categories $D^*_{{\rm add} M}({\rm mod}A)$ and
$D^*_{{\rm Add} M}({\rm Mod}A)$ ({\rm [AHV]}).

\vskip5pt

There is a sequence $D^{b}_{\mathcal C}(\mathcal{A})\subseteq
D^{-}_{\mathcal C}(\mathcal{A})\subseteq D_{\mathcal
C}(\mathcal{A})$ of triangulated subcategories; and the composition
$\mathcal{A} \rightarrow D^b_{\mathcal C}(\mathcal{A})$ of the
embedding $\mathcal{A} \rightarrow K^b(\mathcal A)$ and the
localization functor $K^b(\mathcal A)\rightarrow D^b_{\mathcal
C}(\mathcal{A})$, is fully faithful. The following facts are in
[AHV], [GZ], and [CFH].

\begin{lem} \label{rel} \ Let $\mathcal C$ be a full additive subcategory of abelian category $\mathcal A$. Then

$(1)$ \ Let $C^\bullet\in K^-(\mathcal C)$, and \ $f^\bullet:
X^\bullet \rightarrow C^\bullet$ be a
$\mathcal{C}$-quasi-isomorphism. Then there is $g^\bullet: C^\bullet
\rightarrow X^\bullet$ such that $f^\bullet g^\bullet$ is homotopic
to ${\rm Id}_{C^\bullet}$. Thus, if in addition $X^\bullet\in
K^-(\mathcal C)$, then $f^\bullet$ is a homotopy equivalence.

$(2)$  \ Let $C^\bullet\in K^-(\mathcal C)$ and $Y^\bullet$ be an
arbitrary complex. Then
$\operatorname{Hom}_{K(\mathcal{A})}(C^\bullet, Y^\bullet)\cong
\operatorname{Hom}_{D_{\mathcal C}(\mathcal{A})}(C^\bullet,
Y^\bullet)$, via $f^\bullet\mapsto f^\bullet/{\rm Id}_{C^\bullet}$.
Thus $K^*(\mathcal C)$ can be viewed as a triangulated subcategory
of $D^*_{\mathcal C}(\mathcal{A})$ for $*\in\{b, -\}$.

$(3)$ \ Let $K^{-, \mathcal C, b}(\mathcal C)$ be the full
subcategory of $K^{-}(\mathcal C)$ given by
\begin{align*}\{X^\bullet\in K^{-}(\mathcal C)\ | \ \exists \ N\in\Bbb Z \
\mbox{such that}\ {\rm H}^i\operatorname{Hom}_\mathcal A(C,
X^\bullet) = 0,  \ \forall \ i\le N, \ \forall \ C\in\mathcal
C\}.\end{align*} Then  $K^{-, \mathcal C, b}(\mathcal C)$ is thick
in $K^{-}(\mathcal C)$; and if $\mathcal C$ is contravariantly
finite in $\mathcal A$, then for $X^\bullet \in K^b(\mathcal A)$
there is a $\mathcal{C}$-quasi-isomorphism $C_{_{X^\bullet}}
\rightarrow X^\bullet$  with $C_{_{X^\bullet}}\in K^{-, \mathcal C,
b}(\mathcal C)$.

$(4)$ \ If $\mathcal C$ is a contravariantly finite subcategory of
$\mathcal A$, then there is a triangle-equivalence \ $F: K^{-,
\mathcal C, b}(\mathcal C)\cong D^b_{\mathcal C}(\mathcal A)$ fixing
objects of $K^b(\mathcal C)$.
\end{lem}

\subsection {} Let $\mathcal A$ be an abelian category with enough projective objects, and $\mathcal C$ a resolving subcategory.
Denote by $K^b_{ac}(\mathcal C)$ the full subcategory of
$K^-(\mathcal C)$ consisting of those complexes which are homotopy
equivalent to bounded acyclic complexes over $\mathcal C$. It is
clear that $K^b_{ac}(\mathcal C)$ is a triangulated subcategory of
$K^-(\mathcal C)$.

\begin{lem} \label{b}
Let $C^\bullet\in K^{-, \mathcal C, b}(\mathcal C)$. Then $C^\bullet
\in K^b_{ac}(\mathcal C)$ if and only if $C^\bullet$ is acyclic.
\end{lem} \noindent {Proof.} \  Let $C^\bullet = (C^i,
d^i)$ be acyclic. Since $\mathcal C$ is closed under kernels of
epimorphisms, it follows that $\Ima{d^i} \in \mathcal C, \ \forall \
i\in \Bbb Z$. Since $C^\bullet\in K^{-, \mathcal C, b}(\mathcal C)$,
there is an integer $N$ such that ${\rm H}^{n}\Hom_{\mathcal A}(C,
C^\bullet) = 0, \ \forall \ n \le N, \ \forall \ C \in \mathcal C$.
In particular ${\rm H}^{n}\Hom_{\mathcal A}(\Ima{d^{n-1}},
C^\bullet) = 0.$ This implies that the induced epimorphism
$\widetilde{d^{n-1}}: C^{n-1} \rightarrow \Ima d^{n-1}$ splits for
$n \le N$, and hence there is an isomorphism $C^\bullet\cong
C^{'\bullet}$ in $K^-(\mathcal C)$, where $C'^\bullet\in
K^b_{ac}(\mathcal C)$ is the complex \ $\cdots \rightarrow 0
\rightarrow \Ima d^{N-1} \hookrightarrow C^{N} \rightarrow
C^{N+1}\rightarrow \cdots.$ Thus $C^\bullet \in K^b_{ac}(\mathcal
C)$. $\s$

\begin{prop}\label{relativedescription} Let $\mathcal A$ be an abelian category with enough projective objects, and $\mathcal C$ a resolving contravariantly
finite subcategory. Then $K^b_{ac}(\mathcal C)$ is thick in
$K^{-,\mathcal C, b}(\mathcal C)$, and we have a
triangle-equivalence $D^b(\mathcal A)\cong K^{-,\mathcal C,
b}(\mathcal C)/K^b_{ac}(\mathcal C)$ fixing object $C\in
K^b(\mathcal C)$.
\end{prop}
\noindent{\bf Proof.} Lemma \ref{b} implies that $K^b_{ac}(\mathcal
C)$ is a thick subcategory of $K^{-,\mathcal C, b}(\mathcal C)$. Let
$F': K^b_{ac}(\mathcal C)\rightarrow K^b_{ac}(\mathcal
A)/K^b_{\mathcal Cac}(\mathcal A)$ be the composition of the
embedding $K^b_{ac}(\mathcal C)\hookrightarrow K^{b}_{ac}(\mathcal
A)$ and the Verdier functor $K^{b}_{ac}(\mathcal A)\rightarrow
K^b_{ac}(\mathcal A)/K^b_{\mathcal Cac}(\mathcal A)$. We claim that
$F'$ is an equivalence.

\vskip5pt

In fact, since $K^b_{ac}(\mathcal A)$ is a triangulated subcategory
of $K^b(\mathcal A)$,  $K^b_{ac}(\mathcal A)/K^b_{\mathcal
Cac}(\mathcal A)$ is a triangulated subcategory of the $\mathcal
C$-relative derived category $D^b_\mathcal C(\mathcal A): =
K^b(\mathcal A)/K^b_{\mathcal Cac}(\mathcal A)$. By Lemma
\ref{rel}$(2)$ $F'$ is fully faithful. For $X^\bullet\in
K^{b}_{ac}(\mathcal A)$, by Lemma \ref{rel}$(3)$ there is a
$\mathcal C$-quasi-isomorphism $C^\bullet\rightarrow X^\bullet$ with
$C^\bullet\in K^{-, \mathcal C, b}(\mathcal C)$, which is also a
quasi-isomorphism since $\mathcal C\supseteq \mathcal P$. Since
$X^\bullet$ is acyclic, so is $C^\bullet$. By Lemma \ref{b}
$C^\bullet \in K^b_{ac}(\mathcal C)$. Thus $F'$ is dense, since
$X\cong F'(C^{\bullet})$ in $K^b_{ac}(\mathcal A)/K^b_{\mathcal
Cac}(\mathcal A)$. This proves the claim.

\vskip5pt

By construction $F'$ is the restriction of $F$ to $K^b_{ac}(\mathcal
C)$, where $F$ is the triangle-equivalence $K^{-, \mathcal C,
b}(\mathcal C) \cong D^b_{\mathcal C}(\mathcal A): = K^b(\mathcal
A)/K^b_{\mathcal C ac}(\mathcal A)$ in Lemma \ref{rel}$(4)$. Hence
we get a commutative diagram
\[\xymatrix {K^b_{ac}(\mathcal C)\ar[r]\ar[d]_-{F'} & K^{-,
\mathcal C b}(\mathcal C)\ar[d]_-F
\\ K^b_{ac}(\mathcal
A)/K^b_{\mathcal C ac}(\mathcal A)\ar[r] & K^b(\mathcal
A)/K^b_{\mathcal C ac}(\mathcal A)}\] where the horizontal functors
are embeddings.
Thus $F$ induces a triangle-equivalence
\begin{align*}K^{-,\mathcal Cb}(\mathcal C)/K^b_{ac}(\mathcal C)&\cong (K^b(\mathcal A)/K^b_{\mathcal Cac}(\mathcal A))/(K^b_{ac}(\mathcal
A)/K^b_{\mathcal Cac}(\mathcal A))\\ &\cong K^b(\mathcal A)/K^b_{ac}(\mathcal A) =
D^b(\mathcal A)\end{align*} where the second the triangle-equivalence is well-known (see e.g. [V, Corol. 4-3]). $\s$

\subsection{} Let $\mathcal A$ be an abelian category with enough projective
objects, and  $\mathcal C$ a resolving subcategory. An object $I\in
\mathcal C$ is {\it a {\rm(}relative{\rm)} injective object} of
$\mathcal C$, if $\Hom_\mathcal A(-, I)$ sends any exact sequence
$0\rightarrow X_1\rightarrow X_2\rightarrow X_3\rightarrow 0$ with
each $X_i\in\mathcal C$ to an exact sequence. Clearly, $I$ is an
injective object of $\mathcal C$ if and only if $\Ext_\mathcal
A^1(X, I) = 0$ for $X\in\mathcal C$, and also if and only if
$\Ext_\mathcal A^i(X, I) = 0$ for $X\in\mathcal C$ and for $i\ge 1$.
The following fact is similar to the case of $\mathcal C = \mathcal
A$, which is well-known. Since it will be used in the next section,
we include a justification.

\begin{lem} \label{ff} \ Let $\mathcal A$ be an abelian category with enough
projective objects, and $\mathcal C$ a resolving subcategory of
$\mathcal A$.  Then

$(1)$ \ For $P\in K^-(\mathcal P)$ and $C\in K^{-}(\mathcal C)$,
there is a functorial isomorphism $$\Hom_{K^{-}(\mathcal C)}(P,
C)\cong \Hom_{K^{-}(\mathcal C)/K^{b}_{ac}(\mathcal C)}(P, C), \ \
\mbox {given by} \ f \mapsto f/{\rm Id}_P  \ \mbox {{\rm(}the right
fraction{\rm)}}.$$

$(2)$ \  Let $I$ be a bounded complex of injective objects of
$\mathcal C$ and $G\in K^-_{ac}(\mathcal C)$. Then
$\Hom_{K^-(\mathcal A)}(G, I) = 0.$

 $(3)$ \  Let $I$ be a
bounded complex of injective objects of $\mathcal C$ and $C\in
K^-(\mathcal C)$. If $t: I \rightarrow C$ a quasi-isomorphism, then
there is a chain map $s: C\rightarrow I$ such that $st = {\rm Id}_I$
in $K^-(\mathcal A)$.

$(4)$ \ For a bounded complex $I$ of injective objects of $\mathcal
C$ and $C\in K^{-}(\mathcal C)$, there is a functorial isomorphism
$\Hom_{K^{-}(\mathcal C)}(C, I)\cong \Hom_{K^{-}(\mathcal
C)/K^{b}_{ac}(\mathcal C)}(C, I)$, given by $f \mapsto {\rm
Id}_I\backslash f$ ${\rm(}$the left fraction${\rm)}$.
\end{lem}

\noindent{\bf Proof.} $(1)$ \ Assume $f/{\rm Id}_P = 0$. This means
that there is a quasi-isomorphism $t: Z\rightarrow P$ with $ft = 0$.
Thus there is $g: P \rightarrow Z$ such that $t g$ is homotopic to
${\rm Id}_{P}$ (cf. Lemma \ref{rel}$(1)$). Thus $f = f(t g) = (f t)
g = 0$ in $K^-(\mathcal C)$. Also, assume  $f/s\in
\Hom_{K^-(\mathcal C)/K^b_{ac}(\mathcal C)}(P, C)$, where $s:
Z\rightarrow P$ with ${\rm Con}(s)\in K^b_{ac}(\mathcal C)$. So
there is a quasi-isomorphism $g: P \rightarrow Z$ such that $s g$ is
homotopic to ${\rm Id}_{P}$, and hence  $f/s =fg/sg= fg/{\rm Id}_P$.

\vskip5pt

$(2)$ \ We need to show that any chain map $f: G\rightarrow I$ is
null-homotopic. We construct a homotopy $s = (s^i)$ by induction.
Assume that we have constructed $s^i: G^i \rightarrow I^i$ for $i\le
m$, such that $f^{i-1} = d^{i-2}_Is^{i-1} + s^id_G^{i-1}$ for $i\le
m$. Since $G$ is an upper bounded acyclic complex with all $G^i\in
\mathcal C,$ and $\mathcal C$ is closed under the kernels of
epimorphisms, it follows that ${\rm Im} d^j_G\in \mathcal C, \
\forall \ j\in\Bbb Z.$ Since $(f^m-d_I^{m-1}s^m)d_G^{m-1} = 0$, it
follows that $f^m-d_I^{m-1}s^m$ factors through ${\rm
Coker}d_G^{m-1} = {\rm Im}d^m_G.$ Since $I^m$ is an injective object
of $\mathcal C$, it follows that there is $s^{m+1}:
G^{m+1}\longrightarrow I^m$ such that $f^m-d_I^{m-1}s^m =
s^{m+1}d^m_G.$

\vskip5pt

$(3)$ \ By $(2)$   $\Hom_{K^-(\mathcal A)}({\rm Con} (t), I) = 0 =
\Hom({\rm Con}(t)[-1], I).$ By the distinguished triangle $I
\stackrel{t}\rightarrow C\rightarrow {\rm Con} (t) \rightarrow I[1]$
we see $\Hom(C, I)\cong \Hom(I, I)$, and hence the assertion.

\vskip5pt

 $(4)$ \ The proof is dual to the one of
$(1)$, by using $(3)$. $\s$

\section{\bf Main results}
\begin{thm}\label{mainresult} $(1)$ \ Let $A$ be an Artin algebra
of \ ${\rm gl.dim}A = \infty$, and  $T\in {\rm mod}A$ with \ ${\rm
inj.dim} T< \infty$. If there is a module $M\in{\rm mod}A$ such that
\ $^\perp T = {\rm add} M$ {\rm(}resp. $^{\perp_{\rm big}} ({\rm
Add}T) = {\rm Add} M${\rm)}, then $D^b({\rm mod}B)$ {\rm(}resp.
$D^b({\rm Mod}B)${\rm)} is a categorical resolution of $D^b({\rm
mod}A)$ {\rm(}resp. $D^b({\rm Mod}A)${\rm)}, where $B: = \End M$.

\vskip5pt

$(2)$ \ Let $A$ be a {\rm CM}-finite Gorenstein algebra of ${\rm
gl.dim}A = \infty$, and $B$ its relative Auslander algebra. Then
$D^b({\rm mod}B)$ {\rm(}resp. $D^b({\rm Mod}B)${\rm)} is a weakly
crepant categorical resolution of $D^b({\rm mod}A)$ {\rm(}resp.
$D^b({\rm Mod}A)${\rm)}.
\end{thm}

Note that in general there are no modules $T'\in A$-Mod such that
$^{\perp_{\rm big}} ({\rm Add}T) = \ ^{\perp_{\rm big}} T'$. If in
addition $A$ in Theorem \ref{mainresult}$(2)$ is a commutative local
ring, then G. J. Leuschke [Le] has observed that it is a
non-commutative crepant resolution [VB].

\subsection {} For a right noetherian ring $R$, by [A1] \ ${\rm
gl.dim}({\rm Mod}R) = {\rm gl.dim} ({\rm mod}R),$ which is denoted
by ${\rm gl.dim}R$. For Artin algebra $A$, ${\rm gl.dim}A$ is just
the maximum of ${\rm proj.dim} S(i)$, where $\{S(1), \cdots, S(n)\}$
is the set of simple $A$-modules, up to isomorphisms.

\begin{lem} \label{yoneda} \ Let $A$ be an Artin algebra, and $M\in
$ {\rm mod}$A$. If $M$ is a generator {\rm(}i.e., $A_A\in {\rm
add}M${\rm)}, then $\Hom_{A}(M, -): {\rm mod}A \rightarrow {\rm
mod}B$ is fully faithful.
\end{lem}

\noindent{\bf Proof.} \ Note that $\Hom_{A}(M, -)$ induces an equivalence between $\add M$ and
$\mathcal{P}({\rm mod}B)$ ([ARS, p.33]). Since $M$ is a generator, for any $X\in$ mod$A$ there is a
surjective $A$-map $M^m \twoheadrightarrow X$ for some positive
integer $m$. This implies that $\Hom_{A}(M, -)$ is faithful.

Let $X, Y \in$ mod $A$ and $f: \Hom(M, X)\rightarrow \Hom(M, Y)$ be
a $B$-map. By right $\add M$-approximations we get exact sequences
$T_1 \xrightarrow{u} T_0 \xrightarrow{\pi} X \rightarrow 0$ and
$T_1^\prime \xrightarrow{u^\prime} T_0^\prime
\xrightarrow{\pi^\prime} Y \rightarrow 0$ with $T_0, T_1, T_0',
T_1'\in\add M$. Applying $\Hom(M, -)$ we get the diagram with exact
rows
  \[ \xymatrix{
  \Hom_{A}(M, T_1) \ar[r]^{(M,u)}  \ar@{-->}[d]^{f_1} &\Hom_{A}(M, T_0) \ar[r]^{(M,\pi)} \ar@{-->}[d]^{f_0} &\Hom_{A}(M, X)
  \ar[r] \ar[d]^f &0\\
\Hom_{A}(M, T_1^\prime) \ar[r]^{(M,u^\prime)} &\Hom_{A}(M,
T_0^\prime)
\ar[r]^{(M,\pi^\prime)}  &\Hom_{A}(M, Y) \ar[r] &0.\\
  }\]
Then  $f$ induces $f_1$ and $f_0$ \ such that the diagram
commutes. Thus  $f_i = \Hom_{A}(M, f_i^\prime)$ for some $f_i^\prime\in\Hom_A(T_i, T_i^\prime)$,  $i = 0,
  1$.  So we get the diagram with commutative left square
  \[ \xymatrix{
  T_1 \ar[r]^u  \ar[d]^{f_1^\prime} &T_0 \ar[r]^\pi  \ar[d]^{f_0^\prime}
  &X
  \ar[r] \ar@{-->}[d]^{f^\prime} &0\\
T_1^\prime \ar[r]^{u^\prime}   &T_0^\prime \ar[r]^{\pi^\prime} &Y
  \ar[r]  &0\\
  }\]
and hence there is $f^\prime\in\Hom_A(X, Y)$ such that the diagram commutes. Then one easily deduces that $f = \Hom_{A}(M,
f^\prime)$, i.e., $\Hom_{A}(M, -)$ is full. $\s$

\begin{lem} \label{ABlemma} {\rm (Auslander-Bridger Lemma, [AB, 3.12])} \ Let $\mathcal A$ be an abelian category with enough
projective objects, $\mathcal X$ a resolving subcategory of
$\mathcal A$. For exact sequences
$$0\rightarrow X_n\rightarrow \cdots
\rightarrow X_0 \rightarrow Z\rightarrow 0 \ \ \ \mbox{and}  \ \ \
0\rightarrow Y_n\rightarrow \cdots \rightarrow Y_0 \rightarrow
Z\rightarrow 0$$ in $\mathcal A$ with $X_i\in\mathcal X$ and
$Y_i\in\mathcal X$ for $0\le i\le n-1,$ then $X_n\in\mathcal X$ if
and only if $Y_n\in\mathcal X$.\end{lem}

\begin{prop}\label{global} Let $A$ be an Artin algebra, $T$ and $M$ modules in
${\rm mod}A$ with  \ $^\perp T = {\rm add} M.$ Put $B: = \End M$.
Then for each integer $r\ge 2$,  ${\rm gl.dim} B\le r$ if and only
if ${\rm inj.dim} T\le r$.\end{prop} \noindent{\bf Proof.} \ Assume
that ${\rm gl.dim}B\le r$. Take a right ${\rm add}M$-approximation
$f_0: M_0\rightarrow X$ of $X\in$ mod$A$. Since $M$ is a generator,
$f_0$ is surjective. Again taking a right ${\rm add}M$-approximation
$M_1\rightarrow {\rm Ker}f_0$, and repeating the process we get an
exact sequence
$$M_{r-1} \stackrel{f_{r-1}}\longrightarrow M_{r-2} \rightarrow \cdots
\rightarrow M_0\stackrel {f_0}\longrightarrow  X\rightarrow 0$$ with
each $M_i\in {\rm add}M.$ Put $K: = {\rm Ker}f_{r-1}$. By
construction we get an exact sequence
$$0\rightarrow \Hom(M, K)\rightarrow \Hom(M,
M_{r-1})\rightarrow \cdots \rightarrow \Hom(M, M_0) \rightarrow
\Hom(M, X) \rightarrow 0.$$ Since by assumption ${\rm proj.dim}
_B\Hom(M, X)\le r$, by Auslander-Bridger Lemma $\Hom_A(M, K)$ is a
projective $B$-module. Thus there is a $B$-isomorphism $s: {\rm
Hom}_A(M, K)\rightarrow {\rm Hom}_A(M, M')$ with $M'\in {\rm add}M$.
By Lemma \ref{yoneda} there are $f: K\rightarrow M'$ and $g:
M'\rightarrow K$ such that $s = {\rm Hom}_A(M, f)$ and $s^{-1} =
{\rm Hom}_A(M, g).$ Since $\Hom_A(M, -)$ is faithful, one easily
deduce that $K\in{\rm add}M$, and hence we have an exact sequence
$$0\rightarrow M_r \rightarrow M_{r-1} \rightarrow \cdots
\rightarrow M_0\rightarrow  X\rightarrow 0$$ with each $M_i\in {\rm
add}M.$ Thus for $i\ge 1$ we have $\Ext^{i+r}_A(X, T) \cong
\Ext_A^{i}(M_r, T) = 0$, since $M\in \ ^\perp T$. This implies ${\rm
inj.dim} T \le r$. (Note that this part holds also for $r\le 1$.)

\vskip5pt

Conversely, assume that ${\rm inj.dim} T \le r$ with $r\ge 2$. Let
$_BY\in B$-mod with projective presentation $\Hom_A(M,
M_1)\stackrel d \longrightarrow \Hom_A(M, M_0) \rightarrow \ _BY
\rightarrow 0$ and $M_i\in {\rm add} M, \ i = 0, 1$. Then there
is an $A$-map $f: M_1\rightarrow M_0$ such that $d = \Hom_A(M,
f)$. Taking a right ${\rm add}M$-approximation $M_2\rightarrow {\rm
Ker}f$ and repeating this process we get an exact sequence
$$0\rightarrow K \rightarrow M_{r-1} \rightarrow \cdots
\rightarrow M_2\rightarrow M_1 \stackrel f\longrightarrow M_0$$ with
each $M_i\in {\rm add}M.$ Thus for $i\ge 1$ we have
$\Ext^i_A(K, T) \cong \Ext_A^{i+r}({\rm Coker}f, T) =
0$, since ${\rm inj.dim} T \le r$. So $K\in \ ^\perp T = {\rm add}
M$. By construction we get an exact sequence $$0\rightarrow \Hom(M,
K)\rightarrow \cdots \rightarrow \Hom(M, M_1) \stackrel
d\longrightarrow \Hom(M, M_0) \rightarrow \ _BY\rightarrow 0.$$ This
is a projective resolution, and hence ${\rm proj.dim} _BY\le r$.
This proves ${\rm gl.dim}B\le r$. $\s$

\vskip5pt If $A$ is representation-finite and $T$ an injective
module, then $B$ is the Auslander algebra of $A$ {\rm([ARS])}. If $A$ is {\rm CM}-finite
Gorenstein algebra and  $T: =A_A$, then Theorem \ref{global} is also
well-known, and $B$ is the relative Auslander algebra of $A$ {\rm
(}{\rm[LZ]}, {\rm [Be2];} also {\rm [Le])}.

\subsection{\bf Proof of
Theorem \ref{mainresult}.} $(1)$ \ By Proposition \ref{global}
$D^b({\rm mod}B)$ is smooth. The equivalence $\Hom_A(M, -): {\rm
add} M \rightarrow \mathcal P({\rm mod}B)$ induces pointwisely a
triangle-equivalence $K^{-, {\rm add} M, b}({\rm add} M)\cong K^{-,
b}(\mathcal P({\rm mod}B))$. Since $D^b({\rm mod}B) \cong K^{-,
b}(\mathcal P({\rm mod}B))$, we have a triangle-equivalence
$$F: D^b({\rm mod}B) \cong K^{-, {\rm add} M, b}({\rm add} M).$$
By ${\rm add} M = \ ^\perp T$ we see that ${\rm add} M$ is
resolving, and hence ${\rm add} M$ is a resolving contravariantly
finite subcategory of $A$-mod. By Proposition
\ref{relativedescription} we have a triangle-equivalence
$$G: D^b({\rm mod}A)\longrightarrow K^{-, {\rm add} M, b}({\rm add} M)/K^{b}_{ac}({\rm add} M)$$
with $GP = P$ for $P\in K^b(\mathcal P({\rm mod}A))$. Thus, we get a
triangle functor
$$\pi_*:= G^{-1}V F:  D^b({\rm mod}B)  \longrightarrow
D^b({\rm mod}A),$$ where $V: K^{-, {\rm add} M, b}({\rm add}
M)\rightarrow K^{-, {\rm add} M, b}({\rm add} M)/K^{b}_{ac}({\rm
add} M)$ is the Verdier functor.

On the other hand,  by Proposition \ref{finitely filtrated}$(1)$
$D^b_{\rm perf}({\rm mod}A) = K^b(\mathcal P({\rm mod}A))$. Thus we
have a triangle functor (where $\sigma$ is the embedding $\sigma:
K^b(\mathcal P({\rm mod}A))\hookrightarrow K^{-, {\rm add} M,
b}({\rm add M})$)
$$\pi^*:= F^{-1}\sigma:  D^b_{\rm perf}({\rm mod}A) \longrightarrow
D^b({\rm mod}B).$$

Since $K^b_{ac}({\rm add} M)$ is thick in $K^{-, {\rm add} M,
b}(\mathcal C)$ (cf. Proposition \ref{relativedescription}), we have
${\rm Ker} V = K^b_{ac}({\rm add} M)$. It follows the commutative
diagram \vskip-10pt \[\xymatrix {D^b(B\mbox{-}{\rm
mod})\ar[d]^-{F}\ar[r]^{\pi_*} & D^b({\rm mod}A)\ar[d]^-G
\\K^{-, {\rm add} M, b}({\rm add} M)\ar[r]^-{V} & K^{-, {\rm add} M, b}({\rm add} M)/K^b_{ac}({\rm
add} M)}\] that ${\rm Ker} \pi_* = F^{-} (K^{b}_{ac}({\rm add} M))$,
and $\pi_*$ induces a triangle-equivalence $D^b(B\mbox{-}{\rm
mod})/{\rm Ker}\pi_*\cong D^b({\rm mod}A).$ (But note that $\pi_*$
itself is not full.)

Notice that $\pi^*$ is left adjoint to $\pi_*$ on $K^b(\mathcal
P({\rm mod}A))$. In fact, for $P\in K^b(\mathcal P(A\mbox{-}{\rm
mod}))$ and $X\in D^{b}({\rm mod}B)$ we have $\Hom_{D^b({\rm
mod}B)}(\pi^*P, X)\cong \Hom_{K^{-, {\rm add} M, b}({\rm add} M)}(P,
FX)$ and
\begin{align*}\Hom_{D^b({\rm mod}A)}(P, \pi_*X) & \cong
\Hom_{K^{-, {\rm add} M, b}({\rm add} M)/K^{b}_{ac}({\rm add}
M)}(GP, VFX)\\ & \cong \Hom_{K^{-, {\rm add} M, b}({\rm add}
M)/K^{b}_{ac}({\rm add} M)}(P, FX)\end{align*} (note that $GP = P$
and $VF X = FX$). So, it suffices to prove that there is a
functorial isomorphism
$$\zeta_{P, FX}: \Hom_{K^{-, {\rm add} M, b}({\rm add} M)}(P, FX)\cong \Hom_{K^{-, {\rm add} M, b}({\rm add} M)/K^{b}_{ac}({\rm add} M)}(P,
FX).$$ This follows from Lemma \ref{ff}$(1)$ by taking $\mathcal C:
= {\rm add} M$.

Finally,  saying that the unit ${\rm Id}_{K^b(\mathcal P({\rm
mod}A))} \rightarrow \pi_*\pi^* = G^{-1} V \sigma$ is a natural
isomorphism of functors amounts to saying that $\zeta_P = \zeta_{P,
P}({\rm Id}_P) = {\rm Id}_P/{\rm Id}_P: P\rightarrow P$ is an
isomorphism in $K^{-, {\rm add} M, b}({\rm add} M)/K^{b}_{ac}({\rm
add} M)$ for each $P\in K^b(\mathcal P({\rm mod}A))$. This trivially
holds.

All together the triple $(D^b({\rm mod}B), \pi_*, \pi^*)$ is a
categorical resolution of $D^b({\rm mod}A).$

\vskip5pt

Now we consider $D^b({\rm Mod}A)$. \ First, the condition
$^{\perp_{\rm big}} ({\rm Add}T) = {\rm Add} M$ implies $^\perp T =
{\rm add}M$. The argument is as follows:
$$^\perp T = \ ^{\perp_{\rm big}} ({\rm Add}T) \cap {\rm
mod}A = {\rm Add} M\cap {\rm mod}A = {\rm add}M,$$ where the first
equality follows from the fact that a finitely generated $A$-module
is a compact object in ${\rm Mod}A$. By Proposition \ref{global}
${\rm gl.dim} B < \infty,$ and thus $D^b({\rm Mod}B)$ is smooth.
Since $M$ is finitely generated, it follows that $\Hom_A(M, -): {\rm
Add} M \longrightarrow \mathcal P({\rm Mod}B)$ is again an
equivalence of categories and that it induces pointwisely a
triangle-equivalence $K^{-, {\rm Add} M, b}({\rm Add} M)\cong K^{-,
b}(\mathcal P({\rm Mod}B))$, and hence we get a triangle-equivalence
$$F: D^b({\rm Mod}B) \cong K^{-, {\rm Add} M, b}({\rm Add} M).$$
By ${\rm Add} M = \ ^{\perp_{\rm big}} ({\rm Add} T)$ we see that
${\rm Add} M$ is a resolving subcategory of Mod$A$. Also, ${\rm Add}
M$ is contravariantly finite in Mod$A$ (cf. 1.4). The rest of the
proof is similar with the case of $D^b({\rm mod}A)$, just replacing
${\rm add} M$ by ${\rm Add} M$, and ${\rm mod}B$ by ${\rm Mod}B$.
\vskip5pt

$(2)$ \  In $(1)$ take $T: = A_A$, and $M$ to be the direct sum of
all the pairwise non-isomorphic finitely generated indecomposable
Gorenstein-projective modules.

\vskip5pt

Since $A$ is CM-finite, we have $M\in {\rm mod}A$ and $\mathcal
{GP}({\rm mod}A) = {\rm add} M$. Since $A$ is Gorenstein, it follows
from [EJ, Corol. 11.5.3] that $\mathcal {GP}({\rm mod}A) = \ ^\perp
(A_A)$. Thus $^\perp T = {\rm add} M$, and hence  by $(1)$ \
$D^b({\rm mod}A)$ has a categorical resolution $(D^b({\rm mod}B),
\pi_*, \pi^*)$, where $\pi_*$ and $\pi^*$ are given in the proof of
$(1)$. It remains to see that $\pi^*$ is right adjoint to $\pi_*$ on
$K^b(\mathcal P(A\mbox{-}{\rm mod}))$. In fact, for $X\in D^{b}({\rm
mod}B)$ and $P\in K^b(\mathcal P(A\mbox{-}{\rm mod}))$, as in the
proof of $(1)$ it suffices to prove that there is a functorial
isomorphism
$$\Hom_{K^{-, {\rm add} M, b}({\rm add} M)}(FX, P)\cong \Hom_{K^{-, {\rm add} M, b}({\rm add M})/K^{b}_{ac}({\rm add} M)}(FX, P).$$
This follows from Lemma \ref{ff}$(4)$ by taking $\mathcal C: = {\rm
add} M = \mathcal {GP}({\rm mod}A)$, since projective modules are
injective objects of  $\mathcal {GP}({\rm mod}A)$.

\vskip5pt

Now we consider $D^b({\rm Mod}A)$. Since $A$ is a {\rm CM}-finite
Gorenstein algebra, by [C1] any Gorenstein-projective $A$-module is
a direct sum of finitely generated indecomposable
Gorenstein-projective modules, and hence $\mathcal {GP}({\rm Mod}A)
= {\rm Add} M$. Since $A$ is Gorenstein, it follows from [EJ, Corol.
11.5.3] (or [Be1, Prop. 3.10]) that  $\mathcal {GP}({\rm Mod}A) = \
^{\perp_{\rm big}} ({\rm Add} A_A)$. Thus $^{\perp_{\rm big}} ({\rm
Add} T) = {\rm Add} M,$ and hence by $(1)$ \ $D^b({\rm Mod}B)$ is a
categorical resolution of $D^b({\rm Mod}A)$. The similar argument as
for $D^b({\rm mod}A)$ shows that it is weakly crepant. $\s$

\subsection{} Finally we give some
special cases.

\begin{cor} \label{mainresult1} \  Let $A$ be a representation-finite Artin algebra of ${\rm gl.dim}A = \infty$, and
$B$ its Auslander algebra. Then $D^b({\rm mod}B)$ {\rm(}resp.
$D^b({\rm Mod}B)${\rm)} is a categorical resolution of $D^b({\rm
mod}A)$ {\rm(}resp. $D^b({\rm Mod}A)${\rm)}.
\end{cor}
\noindent{\bf Proof.} Choose $T$ to be an injective module in
mod$A$, and $M$ to be the direct sum of all the pairwise
non-isomorphic finitely generated indecomposable modules.

By Theorem \ref{mainresult}$(1)$ we get the assertion for $D^b({\rm
mod}A)$.

Since $A$ is representation-finite, by [A2, Corol. 4.8] any
$A$-module is a direct sum of finitely generated indecomposable
modules. It follows that $^{\perp_{\rm big}} T = {\rm Mod}A = {\rm
Add} M$, and then the assertion for $D^b({\rm Mod}A)$ follows from
Theorem \ref{mainresult}$(1)$. $\s$

\vskip5pt

A module $T\in$ mod$A$ is {\it cotilting} ([AR]), if ${\rm inj.dim}
T\le 1,$  ${\rm Ext}^1_A(T, T) = 0$, and there is an exact sequence
$0\rightarrow T_0 \rightarrow T_1\rightarrow D(A_A) \rightarrow 0$
with $T_i\in {\rm add} T$, $i =0, 1$.  A module $X\in$ mod$A$ {\it
is cogenerated by} $T$, if $X$ can be embedded as an $A$-module into
a finite direct sum of copies of $T$. Then $X$ is cogenerated by a
cotilting module $T$ if and only if $X\in \ ^\perp T$ ([HR]). By
Theorem \ref{mainresult}$(1)$ we get

\begin{cor} \label{mainresult2} \  Let $A$ be an Artin algebra of ${\rm gl.dim}A = \infty$. If $A$ has a cotilting module $T$ such that
there are only finitely many pairwise non-isomorphic indecomposable
$A$-modules which are cogenerated by $T$. Then $D^b(A\mbox{-}{\rm
mod})$ admits a categorical resolution.
\end{cor}

\end{document}